    \documentclass[11pt,a4paper,oneside]{article}

    \usepackage{amsmath, amsthm, amsfonts, amssymb}
    \usepackage[all]{xy}
    \usepackage{graphics}

	\swapnumbers

    \theoremstyle{plain}
        \newtheorem{theorem}{Theorem}
        \newtheorem{proposition}[theorem]{Proposition}
        \newtheorem{corollary}[theorem]{Corollary}
        \newtheorem{lemma}[theorem]{Lemma}
    \theoremstyle{definition}
        \newtheorem{definition}[theorem]{Definition}
        \newtheorem{remark}[theorem]{Remark}
        \newtheorem{example}[theorem]{Example}

             \def\cat{\mathcal C at}
    \def\top{\mathcal T op}         \def\poset{\mathcal PoSet}
    \def\sset{s\mathcal S et}

    \def\N{\mathbb{N}}
    \def\eps{\varepsilon}

    \def\then{\Rightarrow}                  \def\ls{\leqslant}
                           
    \def\xto#1{\xrightarrow{#1}}

    \def\wt#1{\widetilde{#1}}
    \def\-#1{\overline{#1}}

    \def\title#1{\noindent{\bf\LARGE{#1}} \bigskip \thispagestyle{plain}}
    \def\author#1{\noindent{\sc #1}\smallskip}
    \def\address#1{\noindent #1}

\begin{document}

\title{On the subdivision of small categories}

\author{Matias L. del Hoyo}

\address{Departamento de Matem\'atica\\FCEyN, Universidad de Buenos Aires\\Buenos Aires, Argentina.}

\long\def\symbolfootnote[#1]#2{\begingroup%
\def\thefootnote{\fnsymbol{footnote}}\footnote[#1]{#2}\endgroup}
\symbolfootnote[0]{ {\em 2000 Mathematics Subject Classification:}
18G55; 
55U35; 
18B35. 
}
\symbolfootnote[0]{{\em Key words:} Subdivision; Classifying space; Homotopy category; Posets.}

\begin{abstract}
We present an intrinsic and concrete development of the subdivision of small categories, give some simple examples and derive its fundamental properties. As an application, we deduce an alternative way to compare the homotopy categories of spaces and small categories, by using partially ordered sets. This yields a new conceptual proof to the well-known fact that these two homotopy categories are equivalent.
\end{abstract}

\section*{Introduction}

We present the subdivision of categories from the homotopy point of view, and illustrate this with some simple examples. This subdivision is not new, it has already appeared in some works by Anderson, Dwyer and Kan \cite{anderson,dk}. Here, we derive the basic properties of the subdivision functor $C\mapsto Sd(C)$ from classical results on homotopy of categories, such as the famous Quillen's Theorem A. This way we obtain an intrinsic and geometric-style development of the theory.

\smallskip

Among the fundamental properties of the subdivision of categories, we emphasize theorems \ref{C'' is a poset} and \ref{we}. The first one asserts that
any category becomes a poset after applying the functor $Sd$ twice,
and the second relates the classifying spaces of a category and its subdivision by a homotopy equivalence. These results suggest that the homotopy type of the classifying space of any small category can be modelled by a poset, and therefore that the homotopy categories of small categories and posets are equivalent. This is proved in theorem \ref{fundamental result}.

\smallskip

Finally, we use some results of McCord \cite{mccord} to relate the homotopy categories of posets and topological spaces. Combining these two equivalences we obtain the equivalence of categories
$$Ho(\top)\cong Ho(\cat),$$
which might be thought as a {\em categorical description of topological spaces}. This and the {\em combinatorial description of topological spaces} \cite{gz} are related by Quillen's theorem which asserts that the nerve functor is an equivalence at the homotopy level \cite{illusie}.

\smallskip

I would like to thank my advisor Gabriel Minian for his several suggestions and comments concerning the material of this paper. I also would like to thank Manuel Ladra for his kindness during my days in Santiago de Compostela in which I corrected this work.

%

\section{Preliminaries}


\subsection{Homotopy categories}

If $M$ is a category and $W$ is a family of arrows in $M$, there exists (eventually expanding the base universe) a category $M[W^{-1}]$, called {\em localization of $M$ by $W$}, and a functor $p:M\to M[W^{-1}]$, called {\em localization functor}, which makes invertible all the arrows in $W$ and which is universal for this property.

The category $M[W^{-1}]$ has the same objects than $M$ and its arrows can be expressed as classes of paths involving arrows of $M$ and formal inverses of arrows of $W$ \cite{gz}. By using this description of the localization category, it is easy to prove the following result (cf. \cite{rational}).

\begin{lemma}\label{nt}
Let $p:M\to M[W^{-1}]$ be a categorical localization. Then $p$ induces a bijection
$$Hom(F,G)\xto{\sim}Hom(Fp,Gp)$$
for every pair of functors $F,G:M[W^{-1}]\to N$.
\end{lemma}

When $M$ is a category endowed with homotopical notions (e.g. model categories, simplicial categories, categories with cylinders) and when $W$ is the class of weak equivalences of $M$, the localization of $M$ by $W$ is usually called the {\em homotopy category} of $M$, and is written by $Ho(M)$. The paradigmatic example is that of topological spaces and weak homotopy equivalences. We recall its definition.

\begin{definition}
The homotopy category $Ho(\top)$ is the localization of $\top$, the category of topological spaces, by the family of weak equivalences. Thus, $Ho(\top)=\top[W^{-1}]$ with
$$W=\{f:X\to Y\,|\, f_*:\pi_n(X,x)\to\pi_n(Y,f(x)) \text{ is an isomorphism }\forall n\forall x\}.$$
\end{definition}

\subsection{Homotopical notions in $\cat$}

The other example we are going to consider is that of small categories. The category of small categories $\cat$ is endowed with homotopical notions that one can lift from $\top$ by using the classifying space functor $B:\cat\to\top$ \cite{segal}. We briefly recall from \cite{quillen,segal} some definitions and results concerning this functor.

\smallskip

The category $\Delta$ is that whose objects are the finite ordinals $[q]=\{0<1<...<q\}$ and whose arrows are the order preserving maps.
We use the following standard notation: for $i=0,...,q$ let $s_i:[q+1]\to[q]$ be the surjection which takes twice the value $i$, and let $d_i:[q-1]\to [q]$ be the injection whose image does not contain the value $i$.

If $C$ is a small category, its {\em nerve} $NC$ is the simplicial set whose $q$-simplices are the chains
$$X=(X_0\to X_1\to...\to X_q)$$
of $q$ composable arrows of $C$. Formally, a $q$-simplex $X$ is a functor $[q]\to C$ where $[q]$ is viewed as a category in the canonical way. Faces and degeneracies of $NC$ are given by composing adjacent arrows (or deleting the first or the last arrow) and inserting identities, respectively.

The {\em classifying space} $BC$ of a category $C$ is the geometric realization of its nerve, namely $BC=|NC|$.
A functor $f:B\to C$ in $\cat$ is said to be a {\em weak equivalence} if $Bf$ is a homotopy equivalence in $\top$, and a small category $C$ is said to be {\em contractible} if $BC$ is so.

There is a homeomorphism $B(C\times D)\equiv BC\times BD$ when, for instance, $NC$ or $ND$ has only finite non-degenerate simplices. In particular, denoting $I=[1]$, one has that a functor $C\times I\to D$ induces a continuos map $BC\times BI\to BD$. Thus, it follows that a natural transformation $f\cong g$ induces a homotopy $Bf\cong Bg$. Some simple and useful applications of it are the following.

\begin{lemma}\label{adjoint}
If a functor admits an adjoint, then it is a weak equivalence.
\end{lemma}

\begin{proof}
It suffices to consider the homotopies arising from the unit and the counit of the adjunction.
\end{proof}

\begin{lemma}\label{contractible}
If a category has initial or final object, then it is contractible.
\end{lemma}
\begin{proof}
In these cases the functor $C\to\ast$ admits an adjoint, $\ast$ being the one-arrow category.
\end{proof}


\begin{lemma}\label{retraction}
Let $i:A\to B$ a fully faithful inclusion. If there is a functor $r:B\to A$ and a natural transformation $id_B\then ir$, then $i$ is a weak equivalence.
\end{lemma}

\begin{proof}
The natural transformation $id_B\then ir$ gives rise to another one $i\then iri$, and since $i$ is fully faithful, this is the same than a natural transformation $id_A\then ri$. The result now follows from the fact that a natural transformation induces a homotopy.
\end{proof}

We complete this review by recalling the definition of $Ho(\cat)$.

\begin{definition}
The homotopy category $Ho(\cat)$ is the localization of $\cat$ by the family of weak equivalences, that is, $Ho(\cat)=\cat[W^{-1}]$ with
$$W=\{f:C\to D|Bf:BC\to BD\text{ is a homotopy equivalence}\}.$$
\end{definition}

\begin{remark}
$\cat$ admits a different homotopy structure than the one used here (cf. \cite{minian}). The functors which become homotopy equivalences after taking the classifying space functor are sometimes called {\em topological weak equivalences} to avoid confusions.
\end{remark}

\subsection{Quillen's Theorem A}

Quillen's Theorem A provides a criteria to recognize when a functor is a weak equivalence. We fix some notations and recall it from {\rm \cite[\S 1]{quillen}}.

\smallskip

If $f:C\to D$ is a funtor and if $T$ is an object of $D$, then the {\em fiber} $f^{-1}T$ of $f$ over $T$ is the subcategory of $C$ whose objects and arrows are those which $f$ carries into $T$ and $id_T$, respectively. The {\em left fiber} $f/T$ of $f$ over $T$ is the category of pairs $(X,u)$ with $X$ an object of $C$ and $u:fX\to T$, where an arrow between pairs $(X,u)\to(X',u')$ is a map $v:X\to X'$ in $C$ such that $u'f(v)=u$. The {\em right fiber} $Y/T$ is defined dually. By an abuse of notation, we shall write $C_T$, $C/T$ and $T/C$ for the fiber, left fiber and right fiber, respectively.

\begin{theorem}[Quillen's Theorem A]\label{thm a}
The functor $f:C\to D$ is a weak equivalence if it satisfies either $(i)$ $C/T$ is contractible for every object $T$ of $D$, or $(ii)$ $T/C$ is contractible for every object $T$ of $D$.
\end{theorem}

Given $f:C\to D$, an arrow $X\xto{u} Y$ of $C$ is said to be {\em cocartesian} if every arrow $X\xto{v} Y'$ such that $f(u)=f(v)$ factors as $\wt v\circ u$ with $f(\wt v)=id_{f(Y)}$ in a unique way.
$$\xymatrix@R=10pt{
 & Y'  \\
X \ar[ur]^{\forall v}\ar[r]^u & Y \ar[u]_{\exists!\,\wt v}\\
f(X)\ar[r]^{f(u)} & f(Y)}$$

The functor $f$ is a {\em pre-cofibration} if for each arrow $S\xto{\phi}T$ of $D$ and for each $X\in C_S$ there is a cocartesian arrow $X\xto{u}Y$ over $\phi$. The functor $f$ is a {\em cofibration} if it is a pre-cofibration and also cocartesian arrows are closed under compositions.

{\em Cartesian arrows} are defined dually, as well as {\em pre-fibrations} and {\em fibrations}.

When $f:C\to D$ is a pre-cofibration, the inclusion $C_T\to C/T$ admits a right adjoint, called {\em base-change}, that push-forward an object $(X,\phi)$ along a cocartesian arrow $X\to Y$ over $\phi$. This remark and its dual, combined with lemma \ref{adjoint}, yield the following corollary.

\begin{corollary}\label{coro}
Let $f:C\to D$ be a functor which is either a pre-fibration or a pre-cofibration. If $C_T$ is contractible for every object $T$ of $D$, then $f$ is a weak equivalence.
\end{corollary}

\section{Subdivision of categories}

\subsection{The construction of $Sd(C)$}

Let $C$ be a small category. With $\Delta/C$ we mean the left fiber over $C$ of the embedding $\Delta\to\cat$. It has the simplices of $NC$ as objects, and given $X$ and $Y$ simplices of dimensions $q$ and $p$, a morphism $(Y,\xi,X):X\to Y$ in $\Delta/C$ consists of an order preserving map $\xi:[q]\to [p]$ such that $Y\circ\xi=X$.
We write $\xi_*$ instead of $(Y,\xi,X)$ when there is no place to confusion.

\begin{remark}
Note that if there is a map $X\to Y$ in $\Delta/C$, then the
sequence $X_0\to X_1\to...\to X_q$ is obtained from $Y_0\to
Y_1\to...\to Y_p$ by composing some arrows and inserting some
identities ($X$ is a degeneration of a face of $Y$).
\end{remark}

Let $X\in NC_q$, and let $s:[q+1]\to[q]$ be a surjection. If
$d,d':[q]\to[q+1]$ are the two right inverses of $s$, then we say
that $d_*,d'_*:X\to Xs$ are {\em elementary equivalent}, and
we write $d_*\approx d'_*$. Note that $\approx$ is reflexive and symmetric. We define $\sim$ as the minor
equivalence relation on the arrows of $\Delta/C$ which is
compatible with the composition and satisfies $\xi_*\approx
\xi'_*\then \xi\sim \xi'_*$. We say that $\xi_*$ and $\xi'_*$ are
{\em equivalent} if $\xi_*\sim\xi'_*$. With $[\Delta/C]$ we denote
the quotient category with the same objects than $\Delta/C$ and
arrows the classes under $\sim$.


\begin{definition}
The {\em subdivision of $C$}, denoted $Sd(C)$, is the full subcategory of $[\Delta/C]$ whose objects are the non-degenerate simplices of $NC$.
\end{definition}

We describe the situation with the following diagram,
$$\xymatrix{ & \Delta/C \ar[d] \\ Sd(C) \ar[r] & [\Delta/C] }$$
where $Sd(C)\to\Delta/C$ is just the inclusion and $\Delta/C\to[\Delta/C]$ is the functor which maps an object to itself and an arrow $\xi_*:X\to Y$ to its class $[\xi_*]$ under $\sim$.

\begin{remark}
Notice that this is not the subdivision given in \cite[III-10.1]{adjoint} or \cite[IX-5]{maclane}. Indeed, our construction is equivalent to that of \cite[\S 2]{anderson}, as it can be deduced from lemma \ref{surj2}. Our definition describes completely the arrows of the subdivision category as homotopy-like equivalences of maps, where a degeneration of a simplex plays the role of a cylinder of it.
\end{remark}

\begin{remark}
This subdivision gives rise to a functor $\cat\to\cat$ which equals the composition $c\circ sd\circ N$, where $sd$ denotes Kan's subdivision of simplicial sets \cite{kan}. However, we believe that the intrinsic definition that we present here might be of interest, as it clarifies some aspect of subdivision of categories.
\end{remark}

\begin{example}
If $C$ is the category $\xymatrix@1@C=15pt{0 \ar@<0.7ex>[r]^a \ar@<-0.3ex>[r]_b & 1}$, then the full subcategory of $\Delta/C$ generated by the non-degenerate objects is
$$\xymatrix{ 0 \ar[r] \ar[rd] & a \\ 1 \ar[r] \ar[ru] & b }$$
where $a$ and $b$ denote the non-trivial arrows of $C$. This is also the subdivision $Sd(C)$, since no indentification is possible. The classifying space $B(Sd(C))$ is the 1-sphere $S^1$.
\end{example}

\begin{example}
If $C$ is the two-object simply connected groupoid $0 \rightleftarrows 1$, then $NC$ has two non-degenerate simplices on each dimension $q$, say $0101...$ and $1010...$. If $q<p$, then there are several arrows in $\Delta/C$ between a $q$-simplex and a $p$-simplex, but is not hard to see that any two of them are equivalent. Hence, it follows that $Sd(C)$ is the poset
$$\xymatrix{ 0 \ar[r] \ar[rd] & 01 \ar[r] \ar[rd] & 010 \ar[r] \ar[rd] & \ar@{}[rd]|{...} & \\ 1 \ar[r] \ar[ru] & 10 \ar[r] \ar[ru] & 101 \ar[r] \ar[ru] & & }$$
Notice that $Sd(C)$ is the colimit of its subcategories $Sd(C)_{\leqslant n}$ formed by the simplices of dimension $\leqslant n$. Since $B(Sd(C)_{\leqslant n})=S^n$ and since $B$ commutes with directed colimits, it follows that $B(Sd(C))$ is homeomorphic to the infinite dimensional sphere $S^{\infty}$.
\end{example}


\subsection{Some fundamental properties}

If $X$ is an object of $\Delta/C$, we denote by $q_X$ its dimension as a simplex of $NC$.

\begin{definition}
A map $\xi_*:X\to Y$ in $\Delta/C$ is a {\em surjection} if $\xi:[q_X]\to[q_Y]$ is so.
\end{definition}

Note that if there is a surjection $X\to Y$, then $X$ is a degeneration of $Y$.

\begin{lemma}\label{surj1}
A surjection $\xi_*:X \to Y$ in $\Delta/C$ induces an isomorphism $[\xi_*]:X\to Y$ in $[\Delta/C]$.
\end{lemma}
\begin{proof}
It is sufficient to consider the case $q_X=q_Y+1$, for any surjection can be expressed as a composition of some of the $s_i$.
Thus, suppose that $\xi=s_i:[q+1]\to [q]$, where $q=q_Y$. Then $X=Ys_i$ and the maps $(d_{i+1})_*,(d_{i+2})_*:X\to Xs_{i+1}$ are elementary equivalent.
From the simplicial identities
it follows that
$$(s_i)_*(d_{i+1})_* = id:Y\to Y $$
and that
$$(d_{i+1})_*(s_i)_* = (s_i)_*(d_{i+2})_*\sim(s_i)_*(d_{i+1})_*=id:X\to X.$$
Hence, $(s_i)_*:X\to Y$ and $(d_{i+1})_*:Y\to X$ are inverses modulo equivalences.
\end{proof}

\begin{lemma}\label{surj2}
If a functor $\Delta/C\to D$ carries surjections into isomorphisms, then it factors as $\Delta/C\to[\Delta/C]\to D$ in a unique way. Thus, $[\Delta/C]$ is the localization of $\Delta/C$ by the surjections.
\end{lemma}
\begin{proof}
If it exists, the factorization is unique because $\Delta/C\to[\Delta/C]$ is surjective on objects and on arrows. Let $f:\Delta/C\to D$ be a functor which carries surjections into isomorphisms. If $d_*\approx d'_*:X\to Xs$ are elementary equivalent maps and $s_*:Xs\to X$ is their left inverse, then $f(d_*)=f(s_*)^{-1}=f(d'_*)$. Thus, the relation $\xi_*\sim_f \xi'_* \iff f(\xi_*)=f(\xi'_*)$ is compatible with the composition and satisfies $\xi_*\approx\xi'_*\then \xi_*\sim_f \xi'_*$. Therefore, $\xi_*\sim\xi'_*\then \xi_*\sim_f\xi'_*$ and $f$ factors through $[\Delta/C]$.
\end{proof}

Let $dim:Sd(C)\to \N_0$ be the functor $X\mapsto q_X$ which assigns to each non-degenerate simplex $X$ its dimension.

\begin{lemma}\label{graduada}
If there is a non-identity arrow $X\to Y$ in $Sd(C)$, then $dim(X)<dim(Y)$.
\end{lemma}

\begin{proof}
Let $[i_*]:X\to Y$ be an arrow of $Sd(C)$, with $i:[q_X]\to [q_Y]$ an order preserving map. Then $X=Yi:[q_X]\to C$ and $i$ must be injective because $X$ is a non-degenerate simplex of $NC$. Therefore $q_X\ls q_Y$, and $q_X=q_Y$ if and only if $i=id_{[q_X]}=id_{[q_Y]}$.
\end{proof}

\begin{corollary}\label{no iso}
If $f:X\to Y$ is an isomorphism in $Sd(C)$, then $X=Y$ and $f=id_X$.
\end{corollary}

Following the terminology of \cite[\S 5]{hovey}, we have proved that $dim:Sd(C)\to\N_0$ is a {\em linear extension} of the subdivision category $Sd(C)$, and that the latter is a {\em direct category}. This is not true for $\Delta/C$ nor $[\Delta/C]$, and here lies one reason for our construction.

\begin{theorem}\label{C'' is a poset}
$Sd^2(C)$ is a poset for every category $C$.
\end{theorem}
\begin{proof}
We must show that for every pair $X,Y$ of objects of $Sd^2(C)$, $(i)$ there is at most one arrow $X\to Y$, and $(ii)$ the existence of arrows $X\to Y$ and $Y\to X$ implies $X=Y$.

Assertion $(ii)$ is an immediate corollary of lemma $\ref{graduada}$, so let us prove $(i)$.

An object $X$ of $Sd^2(C)$ is a non-degenerate simplex of $N(Sd(C))$, that is, a chain of non-trivial composable arrows
$$X=(X^0\to X^1\to ...\to X^{q_X})$$
of $Sd(C)$, where $X^i$ is a non-degenerate simplex of $NC$ for each $i$. Note that $dim(X^i)<dim(X^{i+1})$ by lemma \ref{graduada}.

Fix two non-degenerate simplices $X$ and $Y$ of $N(Sd(C))$. We will show that there exists at most one order map $\xi:[q_X]\to[q_Y]$ such that $X=Y\xi$, from where $(i)$ follows.
Suppose that $\xi,\xi'$ are such that $X=Y\xi=Y\xi'$. As we have pointed out, $dim(Y^j)<dim(Y^{j+1})$ for all $j$, so $X^i=Y^{\xi(i)}=Y^{\xi'(i)}$ implies that $\xi(i)=\xi'(i)$ and therefore $\xi=\xi'$.
\end{proof}

\begin{remark}
The same argument of above proves that $Sd(C)$ is a poset for every direct category $C$ in the sense of \cite{hovey}.
\end{remark}

\subsection{Functoriality of the subdivision}

If $X$ is a simplex of $NC$, we might think of $X$ as a sequence of composable arrows, say $X=(f_1,...,f_{q_X})$.
Let $p_X=\#\{j\, |\, f_j\neq id\}$ be the number of non-identity arrows that appear in $X$, and let $r(X)=(f_{i_1},...,f_{i_{p_X}})$ be the sequence obtained from $X$ by deleting the identities. Then $X$ is a degeneration of $r(X)$ -- viewed as simplices of $NC$ -- and $r(X)$ is a non-degenerate $p_X$-simplex. Moreover, $X=r(X)\alpha_X$ with $\alpha_X:[q_X]\to[p_X]$ the surjective order map defined by $\alpha_X(i-1)=\alpha_X(i)\iff f_i=id$.

If $\xi_*:X\to Y$ is an arrow in $\Delta/C$, we define $r(\xi_*)$ as the arrow $r(X)\to r(Y)$ in $Sd(C)$ given by the composition $[(\alpha_Y)_*][\xi_*][(\alpha_X)_*]^{-1}$ in $[\Delta/C]$.
$$\xymatrix{ X \ar[r]^{[\xi_*]} \ar[d]_{[(\alpha_X)_*]} & Y \ar[d]^{[(\alpha_Y)_*]} \\ r(X) \ar@{-->}[r]^{r(\xi_*)} & r(Y) }$$
 Note that $[(\alpha_X)_*]$ is inversible by lemma \ref{surj1}.
 
With these definitions $r=r_C:\Delta/C\to Sd(C)$ is a functor which maps surjection into identities, and by lemma \ref{surj2} it induces a new one $[\Delta/C]\to Sd(C)$, also denoted by $r_C$.

Let $i_C:Sd(C)\to [\Delta/C]$ be the canonical inclusion. Clearly $r_Ci_C=id$, and by lemma \ref{surj1} we have that $\alpha:id\then i_Cr_C$, $X\mapsto [(\alpha_X)_*]$ is a natural isomorphism.
Thus, $i_C:Sd(C)\to [\Delta/C]$ is an equivalence of categories with inverse $r_C:[\Delta/C]\to Sd(C)$.
In particular, $Sd(C)$ is a skeleton of $[\Delta/C]$ as it follows from corollary \ref{no iso}.

\begin{lemma}
The construction $C\mapsto Sd(C)$ is functorial.
\end{lemma}
\begin{proof}
A functor $f:C\to D$ induces a new one $f_*:\Delta/C\to \Delta/D$ by mapping a simplex $X$ to $f\circ X$. This functor clearly sends surjections into surjections. Then, it induces a functor $[f_*]:[\Delta/C]\to[\Delta/D]$, which does not necessarily carry $Sd(C)$ into $Sd(D)$. Thus, we must define $Sd(f)$ as the composition $r_D [f_*] i_C$.

$$\xymatrix{ Sd(C) \ar@{-->}[d]_{Sd(f)} \ar[r]^{i_C} & [\Delta/C] \ar[d]^{[f_*]} \\ Sd(D) & \ar[l]^{r_D} [\Delta/D]}$$

To prove that it is functorial we must verify that $Sd(id)=id$ and that $Sd(g)Sd(f)=Sd(gf)$. The first assertion follows because $Sd(id) = r_C [id_*] i_C = r_C i_C = id_{Sd(C)}$.
About the other, if $f:C\to D$ and $g:D\to E$ then the natural isomorphism $\alpha:id \cong i_D r_D$ induces another one
$$Sd(gf) = r_E [(gf)_*] i_C = r_E [g_*] [f_*] i_C \cong r_E [g_*] i_D r_D [f_*] i_C = Sd(g) Sd(f)$$
of functors $Sd(C)\to Sd(E)$. It follows from corollary \ref{no iso} that the natural isomorphism $Sd(gf)\cong Sd(g)Sd(f)$ must be the identity, and hence we have proved that $Sd(gf)=Sd(g)Sd(f)$ and that $Sd$ is a functor indeed.
\end{proof}

\begin{remark}\label{l}
It follows from proposition \ref{C'' is a poset} that $Sd^2$ lifts to a functor $l:\cat\to\poset$.
$$\xymatrix{ & \poset \ar[d]^j \\ \cat \ar[ru]^l \ar[r]_{Sd^2} & \cat}$$
Here $j$ denotes the canonical inclusion $\poset\to\cat$.  Next section we will show that $l$ is a homotopy inverse for $j$.
\end{remark}

\subsection{Relationship between a category and its subdivision}

Recall the functor $sup:\Delta/C\to C$ \cite[\S 3]{illusie}: Given
an object $X:[q_X]\to C$ of $\Delta/C$, $sup(X)=X_{q_X}$ is the last
object of the sequence $X$. For an arrow $\xi_*:X\to Y$ in
$\Delta/C$, recall that $sup(\xi_*):X_{q_X}\to Y_{q_Y}$ is the
composition of the arrows of $Y$ between $Y_{\xi(q_X)}$ and
$Y_{q_Y}$, namely $sup(\xi_*)=Y(\xi(q_X)\to q_Y)$.

\smallskip

The functor $sup$ maps surjections into identities, since a
surjective map $[p]\to[q]$ preserves final element. It follows
from lemma \ref{surj2} that $sup$ factors through the quotient and
induces a functor $[\Delta/C]\to C$ which will be denoted by
$[sup]$.
$$\xymatrix{ & \Delta/C \ar[d] \ar[dr]^{sup} \\ Sd(C) \ar[r]_{i_C} &  \ar[r]_{[sup]} [\Delta/C]& C}$$

\begin{definition}
The functor $\eps_C:Sd(C)\to C$ is defined as the composition
$[sup]\circ i_C$ of the bottom of the diagram of above.
\end{definition}

\begin{lemma}
The functor $\eps_C$ is natural in $C$. It gives rise to a natural transformation, denoted $\eps:Sd\then id_{\cat}$.
\end{lemma}
\begin{proof}
Given a map $f:C\to D$ in $\cat$, we have the following diagram.
$$\xymatrix{ Sd(C) \ar[d]_{Sd(f)} \ar[r]^{i_C}  & [\Delta/C] \ar[r]^{[sup]} \ar[d]_{[f_*]}& C \ar[d]_f\\
           Sd(D) \ar[r]_{i_D}    & [\Delta/D] \ar[r]_{[sup]}                    & D                 }$$
Clearly the map $sup$ is natural, and it follows from this that
$[sup]$ is also natural. Thus, the right square of above is
commutative. The left square does not commute, but there is a
natural isomorphism $[f_*]i_C\then i_DSd(f)$ which consists of a
surjective map
$$\alpha_{fX}:fX\then i_Dr_D fX$$
for each object $X$ of $Sd(C)$ (see the definition of $\alpha$ in the previous subsection). Finally, as $[sup]$ carries surjections into
identities, the big square commutes and the lemma follows.
\end{proof}



Next we shall prove that the functor $\eps_C:Sd(C)\to C$ is a weak homotopy equivalence. To do that, we first study the left fibers of $\eps_C$. 

\smallskip

Fix some object $T$ of $C$. Let $(X,f)$ be an object of $(\Delta/C)/T$, namely the left fiber of $sup$ over $T$. Thus, $f:S\to T$ is a map in $C$, and $X=(X_0\to X_1\to ...\to X_{q_X-1}\to S)$ is an object of $\Delta/C$ whose top element is $S$. We define $r(X,f)$ as the object of $\Delta/C$ obtained by
extending $X$ with $f$.
$$r(X,f)=(X_0\to X_1\to ...\to X_{q_X-1}\to S \xto{f} T)$$

The assignment $(X,f)\mapsto r(X,f)$ is functorial: given $(X,f)\xto{\xi_*}(Y,g)$, we define $r(\xi_*):r(X,f)\to r(Y,g)$ as the map of $\Delta/C$ induced by the order map
$$[q_X+1]\to[q_Y+1] \qquad j\mapsto \xi(j) \ (0\ls j\ls q_X),\ q_X+1\mapsto q_Y+1.$$
 This way we have a functor $r:(\Delta/C)/T\to (\Delta/C)_T$ into the fiber, which is some kind of retraction for the fully faithful canonical map $i:(\Delta/C)_T\to (\Delta/C)/T$.
Indeed, given $(X,f)$ in $(\Delta/C)/T$, there is a natural map $d_*:X\to r(X,f)$ induced by the injection $d=d_{q_X+1}:[q_X]\to[q_X+1]$. Clearly, $sup(d_*)=f$ and $d_*$ is also a map in $(\Delta/C)/T$, hence we have a natural transformation $id\then ir:(\Delta/C)/T\to(\Delta/C)/T$.

\begin{lemma}
The inclusion $(\Delta/C)_T\to (\Delta/C)/T$ is a weak equivalence.
\end{lemma}
\begin{proof}
Follows from lemma \ref{retraction} and the paragraph of above.
\end{proof}

\begin{remark}\label{to be or not to be}
Note that the map $d_*:X\to r(X,f)$ is not a cocartesian arrow. As an example, consider an arrow $f:X\to T$ in $C$, and let $(X,f)$ be the correspondent zero-dimension object of $(\Delta/C)/T$. Then the two maps
$$[(d_1)_*],[(d_2)_*]:r(X,f)=(X\xto f T)\to(X\xto f T\xto{id} T)$$
are different ways to factor $(X)\to(X\xto f T\xto{id} T)$ through $d_*$.
\end{remark}

Now we shall prove that the functor $r:(\Delta/C)/T\to (\Delta/C)_T$ gives rise to a new one $[r]:[\Delta/C]/T\to [\Delta/C]_T$ between the left fiber and the actual fiber of $[sup]$. In order to do that, we have to show that $r$ carries equivalent maps into equivalent maps.
To prove this, we will need a more explicit description of the relation $\sim$.

\begin{remark}\label{description of sim}
We say that $\xi_*\sim_1\xi'_*$ if there are factorizations $\xi_*=\xi^1_*\xi^2_*...\xi^n_*$ and
$\xi'_*=\xi'^1_*\xi'^2_*...\xi'^n_*$ such that $\xi^i_*\approx\xi'^i_*$ for each $i$. Note that $\sim_1$ is reflexive and symmetric. We call $\sim_2$ to the equivalence relation generated by $\sim_1$.
Thus, $\xi_*\sim_2\xi'_*$ iff there is a sequence
$$\xi_*\sim_1 h_*^1\sim_1h_*^2\sim_1...\sim_1h_*^N\sim_1 \xi'_*.$$
It is easy to see that $\xi_*\sim\xi'_*\iff \xi_*\sim_2\xi'_*$, since $\sim_2$ is an equivalence relation which contains the elementary equivalences and is compatible with the composition (this gives $\then$) and $\sim$ is an equivalence relation which contains $\sim_1$ (this gives $\Leftarrow$).
\end{remark}

\begin{lemma}\label{r preserves sim}
Let $\xi_*,\xi'_*:(X,f)\to (Y,g)$ be maps of $(\Delta/C)/T$ such that $\xi_*\sim\xi'_*$ viewed as maps of $\Delta/C$. Then $r(\xi_*)\sim r(\xi'_*)$.
\end{lemma}

\begin{proof}
First of all, observe that if $\xi_*\approx\xi'_*$, then $r(\xi_*)\approx r(\xi'_*)$.

Secondly, if $\xi_*\sim_1\xi'_*$ then there are factorizations $\xi_*=\xi^1_*\xi^2_*...\xi^n_*$ and
$\xi'_*=\xi'^1_*\xi'^2_*...\xi'^n_*$ such that $\xi^i_*\approx\xi'^i_*$ for each $i$. A priori these are just maps in $\Delta/C$, but since the target of $\xi_*$ and $\xi'_*$ is an object in $(\Delta/C)/T$, then we can think of these maps as arrows in the left fiber. By applying the functor $r$ we obtain factorizations $r(\xi^1_*)r(\xi^2_*)...r(\xi^n_*)$ and $r(\xi'^1_*)r(\xi'^2_*)...r(\xi'^n_*)$ of $r(\xi_*)$ and $r(\xi'_*)$ which together with previous paragraph imply that $r(\xi_*)\sim_1 r(\xi'_*)$.

Finally, if $\xi_*\sim_2\xi'_*$, then  $r(\xi_*)\sim_2 r(\xi'_*)$ by an inductive argument.

The lemma follows from remark \ref{description of sim}.
\end{proof}

\begin{lemma}\label{aux}
The inclusion $[\Delta/C]_T\to [\Delta/C]/T$ is a weak equivalence.
\end{lemma}
\begin{proof}
Follows from lemmas \ref{retraction} and \ref{r preserves sim}.
\end{proof}

The following theorem allow us to consider $Sd(C)$ as an algebraic model for the homotopy type of $BC$, locally simpler than $C$.

\begin{theorem}\label{we}
The functor $\eps_C:Sd(C)\to C$ is a weak equivalence for every $C$.
\end{theorem}

\begin{proof}
The functor $\eps_C$ factors as $[sup]\circ i_C$. Since $i_C$ is an equivalence of categories, it is a weak equivalence (cf. lemma \ref{adjoint}) and we just need to prove that $[sup]$ is a weak equivalence.

We will apply theorem \ref{thm a}, so we need to prove that the left fibers of $[sup]$ are contractible. By lemma \ref{aux} it is sufficient to prove that the fiber $[\Delta/C]_T$ is contractible for each  object $T$ of $C$.

Given $T$, we will prove that $[\Delta/C]_T$ has an initial object and the result will follow from lemma \ref{contractible}. This initial object is $T$, viewed as a 0-simplex of $NC$. If $X$ is any object of $[\Delta/C]_T$, then the $(q_X)$-th inclusion $\alpha:[0]\to[q_X]$ induces a map $[\alpha_*]:T\to X$ in $[\Delta/C]_T$.

If $[\beta_*]:T\to X$ is any other map in $[\Delta/C]_T$, we have to prove that $\alpha_*\sim\beta_*$. Consider the order map $h:[1]\to[q_X]$ given by $h(0)=\beta(0)$ and $h(1)=q_X$.
$$\begin{matrix}\xymatrix@C=50pt{
(T\xto{id}T) \ar@/^/[dr]^{h_*} \ar[d]|{(s_0)_*}& \\
(T) \ar@<2.6ex>[u]^{(d_0)_*} \ar@<-2.2ex>[u]_{(d_1)_*} \ar@<1ex>[r]^(.4){\alpha_*}\ar@<-1ex>[r]_(.4){\beta_*}
& (X_0\to X_1\to...\to X_{q_X})}\end{matrix}\qquad \begin{matrix}h_*(d_0)_*=\alpha_* \\ h_*(d_1)_*=\beta_*\end{matrix}
$$
Then $\alpha=hd_0$ and $\beta=hd_1$, and because $X(\beta(0)\to q_X)=\eps_C([\beta_*])=id_T$ it follows that $Xh=Ts_0$ is a degeneration of $T$, $(d_0)_*\approx(d_1)_*$ and therefore $\alpha_*\sim\beta_*$.
\end{proof}

\begin{corollary}\label{sd preserva we}
The functor $Sd:\cat\to\cat$ preserves weak equivalences.
\end{corollary}
\begin{proof}
If $f:C\to D$ is a weak equivalence in $\cat$, it follows from theorem \ref{we} and the square
$$\xymatrix{Sd(C) \ar[r]^{\eps_C}\ar[d]_{Sd(f)} & C \ar[d]^f \\ Sd(D) \ar[r]_{\eps_D} & D}$$
that $Sd(f)$ is also a weak equivalence.
\end{proof}

\begin{remark}
Given $(X,f)$ an object of $(\Delta/C)/T$, we have seen in remark \ref{to be or not to be} that $d_*:X\to r(X,f)$ is not a cartesian arrow for $sup$. However, $[d_*]:X\to r(X,f)$ is a cocartesian arrow for $[sup]$. To see that, suppose that $[\xi_*]:X\to Y$ is an arrow of $[\Delta/C]$ such that $[sup]([\xi_*])=f$. Then, $\xi_*$ might be consider as an arrow $(X,f)\to(Y|_{[q_Y-1]},Y([q_Y-1]\to[q_Y]))$ in $[\Delta/C]/T$, and $[\xi_*]$ factors as $[r(\xi_*)][d_*]$ (actually, $\xi_*=r(\xi_*)d_*$). To see that this factorization is unique, suppose that another one is given, and use the fact that $r$ preserves equivalences.

It follows that $[\Delta/C]\to C$ is a precofibration, as well as $Sd(C)\to C$. Thus, theorem \ref{we} can be proved by using corollary \ref{coro}. However, $Sd(C)\to C$ is not a cofibration in general, since cocartesian arrows are not closed under composition. This is clear because a cocartesian arrow over a non-identity map must increase the degree in exactly one.
\end{remark}

\section{Application to homotopy theory}

\subsection{Homotopy category of $\poset$}

Despite the homotopy theory of partially ordered sets is largely developed, we could not find a definition for the homotopy category $Ho(\poset)$. We construct it here in a suitable form, compatible with the inclusions $\poset\to\top$ and $\poset\to\top$.
We will use for this purpose some well known facts about $A$-spaces, posets and simplicial complexes.

\smallskip

Recall that an {\em A-space}, or {\em Alexandrov space}, is a topological space in which any arbitrary intersection of open subsets is open. A topological space satisfies the $T_0$ separability axiom if given two points on it, there exists an open subset that contains exactly one of these points. A $T_0A$-space is simply an $A$-space which satisfies the $T_0$ axiom.

There is a well-known correspondence between $T_0A$-spaces and preorders. We recall it briefly.

If $P$ is a poset, let $a(P)$ be the topological space with points the elements of $P$ and with open basis formed by the subsets $\{y|y\ls x\}$, $x\in P$. Clearly, $a(P)$ is a $T_0A$-space.

If $X$ is a $T_0A$-space, let $s(X)$ be the poset with elements the points of $X$ and with the order $x\ls y \iff y\in cl(x)$, where $cl(x)$ denotes the closure of $\{x\}$ in $X$. Note that the relation $\ls$ is antisymmetric because $X$ is $T_0$.

\begin{lemma}
The constructions $P\mapsto a(P)$ and $X\mapsto s(X)$ are functorial, and they define an equivalence of categories between $\poset$ and the full subcategory of $\top$ whose objects are the $T_0A$-spaces.
\end{lemma}

We recall some constructions from \cite{mccord}.
Given $X$ a $T_0A$-space, a simplicial complex $k(X)$ is constructed with vertices the points of $X$ and simplices the finite chains of $s(X)$, namely the sequences of points $(x_0,...,x_q)$ satisfying $x_{i+1}\in cl(x_i)$. The construction $X\mapsto k(X)$ is functorial. Moreover, there is a natural continuous map $f_X:|k(X)|\to X$ defined by $f_X(u)=min(carrier(u))$, where $carrier(u)$ is the unique open simplex containing $u$.

Given a simplicial complex $K$, denote by $S(K)$ its set of simplices ordered by inclusion. Define $x(K)$ as the $T_0A$-space associated to its simplices, namely $x(K)=aS(K)$. The construction $K\mapsto x(K)$ is functorial, since $a$ and $S$ are so. Moreover, since $k(x(K))$ is just the barycentric subdivision of $K$, there is a natural continuous map $f_K:|K|\to x(K)$ defined as the composition of the canonical homeomorphism $|K|\xto{\sim}|kx(K)|$ with the map $f_{x(K)}$.
The following results are due to  McCord \cite{mccord}.

\begin{proposition}\label{teo2}
For every $T_0A$-space $X$ the map $f_X:|k(X)|\to X$ is a weak homotopy equivalence.
\end{proposition}

\begin{proposition}\label{mc}
For every simplicial complex $K$ the map $f_K:|K|\to x(K)$ is a weak homotopy equivalence.
\end{proposition}

Now we are in condition to describe $Ho(\poset)$.
Recall that $j:\poset\to\cat$ is the functor which assigns to each poset $P$ a category $j(P)$ in the usual way. The functor $j$ admits a left adjoint $p:\cat\to\poset$, which assigns to each small category $C$ the poset associated to the preorder defined over the objects of $C$ by the rule
$$X\ls Y \iff \text{ there exists an arrow } X\to Y.$$

The functors $a:\poset\to\top$ and $j:\poset\to\cat$ embed $\poset$ as a full reflective subcategory of $\top$ and $\cat$. Thus, $\poset$ inherits two definitions for weak equivalences by lifting those of $\top$ and $\cat$.
Let $W_a$  be the class of maps $f:P\to Q$ in $\poset$ such that $a(f):a(P)\to a(Q)$ is a weak equivalence in $\top$, and let $W_j$ be the class of maps $f:P\to Q$ in $\poset$ such that $j(f):j(P)\to j(Q)$ is a weak equivalence in $\cat$ or, what is the same, $Bj(f):Bj(P)\to Bj(Q)$ is a weak equivalence in $\top$.

\begin{proposition}
The classes $W_a$ and $W_j$ coincide.
\end{proposition}

\begin{proof}
For each poset $P$ there is a natural homeomorphism $Bj(P)\cong |ka(P)|$ between the classifying space of $j(P)$ and the geometric realization of McCord's construction on $a(P)$.
Given $f:P\to Q$ a map in $\poset$, consider the following commutative diagram.
$$\xymatrix{
Bj(P) \ar[d]_{Bj(f)} \ar[r]^{\sim} & |ka(P)| \ar[r]^{f_{a(P)}} \ar[d]_{|Ka(f)|} &  a(P) \ar[d]_{a(f)}  \\
Bj(Q)                \ar[r]^{\sim} & |ka(Q)| \ar[r]_{f_{a(Q)}}                  &  a(Q)                }$$
Since the maps $f_{a(P)}$ and $f_{a(Q)}$ are weak equivalences in $\top$ (cf. proposition \ref{teo2}), the continuous map $Bj(f)$ is a weak equivalence if and only if $a(f)$ is so.
\end{proof}

\begin{definition}
We say that a map $f:P\to Q$ in $\poset$ is a {\em weak equivalence} if $f\in W_a=W_j$.
We define the {\em homotopy category of $\poset$}, denoted $Ho(\poset)$, as the localization of $\poset$ by the family of weak equivalences.
\end{definition}

\begin{remark}
It is clear that $pj(P)=P$. Unfortunately, the composition $jp$ does not preserve homotopy types -- for instance, a group $G$ is mapped by $jp$ into the one-arrow category.
Similarly, while the composition $sa$ is the identity functor over $\poset$, the other composition $as$ fails at the homotopy level -- for instance, a Hausdorff space $X$ is mapped by $as$ into a discrete space.
\end{remark}

Despite last remark, the functors $j:\poset\to\cat$ and $a:\poset\to\top$ induce equivalences between the homotopy categories. In the next subsection we will construct homotopy inverses to the inclusions $a$ and $j$.

\subsection{Categorical description of $Ho(\top)$}

The functors $a$ and $j$ preserve weak equivalences. Hence, they induce functors $Ho(a)$ and $Ho(j)$ at the homotopy level.
$$\xymatrix{ \cat \ar[d]& \poset \ar[r]^a \ar[d]\ar[l]_j & \top \ar[d] \\
 Ho(\cat) &  Ho(\poset) \ar[r]_{Ho(a)} \ar[l]^{Ho(j)} & Ho(\top)}$$

\begin{theorem}\label{fundamental result}
The functors $Ho(a)$ and $Ho(j)$ are equivalences of categories.
Hence, the categories $Ho(\cat)$ and $Ho(\top)$ are equivalent.
\end{theorem}

This theorem is intimately related with Quillen's theorem asserting that $N$ induces an equivalence of categories at the homotopy level (cf. \cite{illusie}). One can derive one from the other by using the well known equivalence $Ho(\top)\cong Ho(\sset)$.

\begin{proof}
We prove first that $Ho(j)$ is an equivalence of categories. Recall from remark \ref{l} the definition of $l:\cat\to\poset$.
We have seen in corollary \ref{sd preserva we} that $Sd$ preserves weak equivalences. Since $jl=Sd^2$, it is clear that $l$ preserves them too, hence it induces a functor $Ho(l):Ho(\cat)\to Ho(\poset)$. We assert that $l$ is a homotopy inverse to $j$, so we have to prove that there are
natural isomorphisms $Ho(jl)=Ho(j)Ho(l)\cong id_{Ho(\cat)}$ and $Ho(lj)=Ho(l)Ho(j)\cong id_{Ho(\poset)}$.

If we show that there are natural transformations $jl\then id_{\cat}$ and $lj\then id_{\poset}$ which assign to any object a weak equivalence, then by composing with the projections we will obtain natural isomorphisms, which yield another ones $Ho(jl)\cong id_{Ho(\cat)}$ and $Ho(lj)\cong id_{Ho(\poset)}$ by lemma \ref{nt}.

For every category $C$ the composition $\eps_C\eps_{Sd(C)}:jl(C)=Sd^2(C)\to C$ is a weak equivalence by theorem \ref{we}, and clearly it is natural. This gives the natural isomorphism $Ho(jl)\cong id_{Ho(\cat)}$. The other natural isomorphism can be obtained as a restriction of this.



\smallskip

Now we prove that $Ho(a)$ is an equivalence of categories.
We will construct an inverse to $a$ by considering for each topological space $X$ a simplicial complex $K_X$ and a weak equivalence $|K_X|\to X$, which can be done naturally. We define a functor $b:\top\to\poset$ by $b(X)=S(K_X)$ the poset of simplices of the associated complex. To see that $b$ preserves weak equivalences it is sufficient to consider the diagram
$$\xymatrix{ab(X)=aS(K_X) \ar[d] & |K_X|  \ar[l]_(.3){f_{K_X}}\ar[r]^{\sim} \ar[d] & X  \ar[d]\\ ab(Y)=aS(K_Y) & |K_Y| \ar[l]^(.3){f_{K_Y}} \ar[r]_{\sim} & Y }$$
where $f_{K_X}$ and $f_{K_Y}$ are McCord's weak equivalences of proposition \ref{mc}.
Hence $b$ induces a functor $Ho(b):Ho(\top)\to Ho(\poset)$.

By the same argument used above, the natural weak equivalences $|K_X|\to X$ and $|K_X|\to ab(X)$ yield natural isomorphisms at the homotopy level, which compose to give $Ho(a)Ho(b)=Ho(ab)\cong id_{Ho(\top)}$.
The natural isomorphism $Ho(b)Ho(a)=Ho(ba)\cong id_{Ho(\poset)}$ can be obtained as a restriction of the previous one.
\end{proof}

\begin{remark}
By the work of Thomason \cite{thomason} we know that $\cat$ admits a closed model structure, weak equivalences being the ones we work with. By the corrections made by Cisinski \cite{cisinski} over the paper of Thomason, we know that every cofibrant category under this structure is a poset. Thus, the equivalence $Ho(\poset)\xto{\sim} Ho(\cat)$ can be deduced from the composition $$Ho(\cat_c) \xto{\sim} Ho(\poset) \xto{\sim} Ho(\cat),$$ where $\cat_c$ denotes the full subcategory of cofibrant objects.
\end{remark}


\end{document}